\UseRawInputEncoding
\documentclass[11pt]{article}
\usepackage{amsfonts}
\usepackage{color}
\usepackage{amsmath,amssymb,comment}
\usepackage{verbatim}
\newtheorem{theo}{Theorem}[section]
\newtheorem{lem}[theo]{Lemma}

\newtheorem{defi}[theo]{Definition}

\newcommand{\mysection}[1]{\section{#1} \setcounter{equation}{0}}
    \makeatletter
\def\@fnsymbol#1{\ensuremath{\ifcase#1\or *\or \ddagger\or
   \mathsection\or \mathparagraph\or \|\or **\or \dagger\dagger
   \or \ddagger\ddagger \else\@ctrerr\fi}}
    \makeatother
\newcommand{\proof}{{\sc Proof.} \quad}
\newcommand{\proofc}{{\sc Proof} \ }
\newcommand{\be}{\begin{equation} \label}
\newcommand{\ee}{\end{equation}}
\newcommand{\bea}{\begin{eqnarray}\label}
\newcommand{\eea}{\end{eqnarray}}
\newcommand{\bas}{\begin{eqnarray*}}
\newcommand{\eas}{\end{eqnarray*}}
\newcommand{\bit}{\begin{itemize}}
\newcommand{\eit}{\end{itemize}}
\newcommand{\qed}{\hfill$\Box$ \vskip.2cm}
\newcommand{\nn}{\nonumber}
\newcommand{\R}{\mathbb{R}}
\newcommand{\N}{\mathbb{N}}
\newcommand{\pO}{\partial\Omega}

\newcommand{\eps}{\varepsilon}

\newcommand{\wto}{\rightharpoonup}
\newcommand{\wsto}{\stackrel{\star}{\rightharpoonup}}
\newcommand{\hra}{\hookrightarrow}
\newcommand{\io}{\int_\Omega}
\newcommand{\na}{\nabla}
\newcommand{\Del}{\Delta}
\newcommand{\del}{\delta}

\newcommand{\pa}{\partial}
\newcommand{\bom}{\overline{\Omega}}
\newcommand{\Om}{\Omega}

\newcommand{\ov}{\overline}

\newcommand{\hs}{\hspace*}
\newcommand{\sm}{\setminus}

\newcommand{\vp}{\varphi}
\newcommand{\lbal}{\left\{ \begin{array}{l}}
\newcommand{\lball}{\left\{ \begin{array}{ll}}
\newcommand{\ear}{\end{array} \right.}
\newcommand{\ouz}{\ov{u}_0}

\newcommand{\abs}{\\[5pt]}

\newcommand{\tme}{T_{max,\eps}}

\newcommand{\pe}{\phi_\eps}
\newcommand{\one}{{\mathbf{1}}}

\newcommand{\oy}{\ov{y}}

\newcommand{\ueps}{u_\eps}
\newcommand{\veps}{v_\eps}

\newcommand{\heps}{h_\eps}

\newcommand{\geps}{g_\eps}
\newcommand{\yeps}{y_\eps}

\newcommand{\Feps}{{\mathcal{F}}_\eps}
\newcommand{\F}{{\mathcal{F}}}
\begin{document}
\enlargethispage{10mm}
\title{Relaxation in a Keller-Segel-consumption system\\
involving signal-dependent motilities}
\author{
Genglin Li\footnote{1182028@mail.dhu.edu.cn}\\
{\small College of Information and Technology, }\\
{\small Donghua University, Shanghai 201620, P.R.~China}
\and
Michael Winkler\footnote{michael.winkler@math.uni-paderborn.de}\\
{\small Institut f\"ur Mathematik, Universit\"at Paderborn,}\\
{\small 33098 Paderborn, Germany} 
}
\date{}
\maketitle
\begin{abstract}
\noindent 
%
%
%This paper is concerned with 
Two relaxation features of the migration-consumption chemotaxis system involving signal-dependent motilities,
\bas
	\left\{ \begin{array}{l}
	u_t = \Delta \big(u\phi(v)\big), \\[1mm]
	v_t = \Delta v-uv,
	\ear
	\qquad \qquad (\star)
  \eas
are studied in smoothly bounded domains $\Om\subset\R^n$, $n\ge 1$:
% for which the non-degenerate diffusion, namely, $\phi>0$ on $[0,\infty)$, is considered. 
%This study mainly focuses on the construction of global (very) weak solutions to the system
%$(\star)$ as well as on their large time stabilization, in the case when the motility function $\phi$ is assumed to be 
%suitably smooth, particularly the extent to which the relaxation of the regularity hypotheses on initial data of arbitrary size 
%occurs. In this paper i
It is shown that if $\phi\in C^0([0,\infty))$ is positive on $[0,\infty)$, then for any initial data $(u_0,v_0)$ belonging to 
the space $(C^0(\bom))^\star\times L^\infty(\Om)$ an associated no-flux type initial-boundary value problem admits a 
global very weak solution.
Beyond this initial relaxation property, it is seen that
under the additional hypotheses that $\phi\in C^1([0,\infty))$ and $n\le 3$,
each of these solutions stabilizes toward a semi-trivial spatially homogeneous steady state in the large time limit.\abs
By thus applying to irregular and partially even measure-type initial data of arbitrary size, 
this firstly extends previous results on global solvability in ($\star$) which have been restricted to initial data
not only considerably more regular but also suitably small.
Secondly, this reveals a significant difference between the large time behavior in ($\star$) and that in related 
degenerate counterparts involving functions $\phi$ with $\phi(0)=0$, about which, namely, it is known that some
solutions may asymptotically approach nonhomogeneous states.\abs
%
%
%This seemingly provides a considerably wider range of the eligible set of initial data allowing for the construction of global solutions 
%to $(\star)$ than does a recent result obtained for small data with further restrictions on the regularity. \abs
%
%
\end{abstract}
\newpage
\section{Introduction}\label{intro}
Chemotaxis systems accounting for local sensing mechanisms have attracted increased interest in the recent literature.
When viewed in the context of general Keller-Segel type systems (\cite{KS70}) of the form
\be{ks}
	\left\{ \begin{array}{l}	
	u_t = \na\cdot \big( D(u,v)\na u - \chi(u,v)u\na v\big), \\[1mm]
	v_t = \Delta v+K(u,v),
	\end{array} \right.
\ee
setting a corresponding focus amounts to assuming the diffusivity $D$ and the cross-diffusion rate $\chi$ to be linked through
the relations 
\be{dc}
	D(u,v)=\phi(v)
	\qquad \mbox{and} \qquad
	\chi(u,v)=\phi'(v), 
\ee
with some nonnegative function $\phi$ exhibiting
suitable decay at large values of the signal concentration $v$ in order to reflect the local character of sensing
(\cite{fu}, \cite{liu}); typical examples thus include 
\be{typ}
	\phi(\xi)=\frac{1}{(\xi+a)^\alpha}
	\qquad \mbox{or also} \qquad
	\phi(\xi)=e^{-\beta\xi}
\ee
for $\xi\ge 0$, with $\alpha>0, a\ge 0$ and $\beta>0$.
For accordingly obtained versions of (\ref{ks}) that address situations in which the considered signal is produced
by cells, remarkably comprehensive knowledge has been achieved in several respects.
Namely, in the case when $K(u,v)=u-v$, throughout considerably large classes of the key ingredient $\phi$ to (\ref{ks})-(\ref{dc})
the literature meanwhile provides not only far-reaching results on global solvability 
(\cite{fujie_senba}, \cite{ahn_yoon}, 
\cite{fujie_jiang_ACAP2021}, \cite{taowin_M3AS}, \cite{burger}, \cite{desvillettes}; cf.~also
\cite{jiang_laurencot}, \cite{jin_kim_wang}, \cite{wang_wang}, \cite{yifu_wang}, \cite{wenbin_lv_ZAMP}, \cite{wenbin_lv_EECT}
and \cite{wenbin_lv_PROCA}	%and \cite{win_NON} 
for corresponding studies on some close relatives), 
but also on large time asymptotics, and especially on the identification of situations in which
either stabilization toward equilibria can be observed, or, alternatively, infinite-time blow-up occurs 
(\cite{DLTW}, \cite{taowin_M3AS}, \cite{fujie_jiang_CVPDE}, \cite{fujie_senba}, \cite{jin_wang}).
In particular, it has been found that within this class of chemotaxis-production systems,
appropriate choices of the form in (\ref{typ}) lead to a substantial support of spatial structures in the sense either of
large-time singularity formation (\cite{fujie_jiang_CVPDE}, \cite{fujie_senba}, \cite{jin_wang}), or at least 
of heterogeneous long-term aymptotics trivially exhibited by non-constant steady state solutions (\cite{DLTW}).
These observations are quite in line with numerous findings on emergence and stabilization of singular structures
(\cite{herrero_velazquez}, \cite{nagai2001}, \cite{suzuki_book}, \cite{senba_suzuki}, \cite{biler2008}, 
\cite{cieslak_stinner_JDE2015}),		%, \cite{win_JMPA}),
and supplementary also of more subtle destabilization of spatial homogeneity (\cite{win_MATANN}),	%, \cite{ct_stable2}), 
in various further among the chemotaxis-production versions of (\ref{ks}) with constituents more general than in (\ref{dc}).\abs
In contrast to this, the present study now focuses on contexts in which chemotactic motion based on local sensing 
%of the form described by (\ref{ks})-(\ref{dc}) 
is directed by a cue that is consumed by individuals, rather than produced.
By describing signal absorption in the apparently most standard functional form, a resulting chemotaxis-consumption 
version of (\ref{ks})-(\ref{dc}) becomes
\be{00}
	\left\{ \begin{array}{l}	
	u_t = \Del \big(u\phi(v)\big), \\[1mm]
	v_t = \Delta v -uv,
	\end{array} \right.
\ee
and in stark difference to the situations discussed above, 
questions related to the evolution of structures in the above sense seem to have remained widely unaddressed in 
such frameworks.		% of this type.
Indeed, the only result concerned with global solutions to a problem of this form in high-dimensional settings,
as recently achieved in \cite{li_zhao_ZAMP}, relies on a smallness restriction on the initial data to assert 
global classical solvability and large time convergence to homogeneous states in non-degenerate cases in which 
$\phi$ is strictly positive throughout $[0,\infty)$.
Large-data solutions appear to have been constructed only in spatially one- and two-dimensional domains in
\cite{win_sig_dep_mot_cons_2}, where a focus has been on a degenerate version in which
yet $\phi>0$ on $(0,\infty)$ but
$\phi(0)=0$; in such cases, for suitably regular initial data some global smooth solutions have found to exist 
and to approach a steady state $(u_\infty,0)$ in the large time limit, with $u_\infty$ known to be nonconstant whenever
$v|_{t=0}$ is appropriately small (\cite{win_sig_dep_mot_cons_2}).\abs
{\bf Main results.} \quad
One purpose of the present study is to make sure that this latter observation, paralleling similar findings on the existence
of non-constant large-time patterns in related taxis-consumption systems of the form (\ref{ks}) with signal-dependent
motility degeneracies (\cite{win_TRAN2021}),
cannot be made in any non-degenerate version of (\ref{00}).
Beyond asserting this absence of structure support on large time scales, however, a second goal will consist in identifying
a second pattern-countercting feature of (\ref{00}) which will become manifest in a result on instantaneous relaxation
of solutions emanating from considerably irregular initial data.\abs
%
%To make this more precise,
More precisely,
in a smoothly bounded domain $\Om\subset \R^n$, $n\ge 1$, let us consider the initial-boundary value problem
\be{0}
	\left\{ \begin{array}{ll}	
	u_t = \Del \big(u\phi(v)\big),
	\qquad & x\in\Om, \ t>0, \\[1mm]
	v_t = \Delta v-uv,
	\qquad & x\in\Om, \ t>0, \\[1mm]
	\na \big(u\phi(v)\big) \cdot \nu = \na v \cdot \nu=0,
	\qquad & x\in\pO, \ t>0, \\[1mm]
	u(x,0)=u_0(x), \quad v(x,0)=v_0(x),
	\qquad & x\in\Om,
	\end{array} \right.
\ee
and first concentrate on the latter question on initial relaxation
by deriving a result on global solvability under mild requirements on data regularity which inter alia even allow for measure-valued
first components of the initial distributions.
In formulating this and throughout the sequel, 
we let $C^0_{w-\star}([0,\infty);(C^0(\bom))^\star)$ and $C^0_{w-\star}([0,\infty);L^\infty(\Om))$
denote the spaces of functions which are continuous on $[0,\infty)$ with respect to the weak-$\star$ topology in
$(C^0(\bom))^\star$ and $L^\infty(\Om)$, respectively.\abs
Specifically, the first of our main results indeed reveals that 
%unlike in virtually all precedent results on Keller-Segel type
%problems involving irregular initial data, 
even such singular initial settings undergo
an immediate relaxation into globally existing solutions to (\ref{0}) which belong to
$L^2(\Om) \times W^{2,2}(\Om)$ at a.e.~positive time, under the mere assumption that $\phi$ be continuous and positive:
\begin{theo}\label{theo10}
  Let $n\ge 1$ and $\Om\subset\R^n$ be a bounded domain with smooth boundary, and suppose that
  \be{phi}
	\phi\in C^0([0,\infty))
	\mbox{\quad is such that \quad } \phi(\xi)>0 \mbox{ for all } \xi\ge 0,
  \ee
  and that 		%$(u_0,v_0)$ complies with (\ref{init}).
  \be{init}
	u_0\in (C^0(\bom))^\star
	\mbox{ as well as }
	v_0\in L^\infty(\Om)
	\quad \mbox{are nonnegative.}
  \ee
  Then there exist nonnegative functions
  \be{reg}
	\lbal
	u\in C^0_{w-\star}([0,\infty);(C^0(\bom))^\star) \cap L^\infty((0,\infty);L^1(\Om)) \cap L^2_{loc}(\bom\times (0,\infty))
	\qquad \mbox{and} \\[1mm]
	v\in C^0_{w-\star}([0,\infty;L^\infty(\Om)) \cap L^\infty(\Om\times (0,\infty)) \\
	\hs{10mm}
	\cap L^2_{loc}((0,\infty);W^{2,2}(\Om)) \cap L^4_{loc}((0,\infty);W^{1,4}(\Om)) \cap L^\infty_{loc}((0,\infty);W^{1,2}(\Om))
	\ear
  \ee
  which are such that $(u,v)$ forms a global very weak solution of (\ref{0}) in the sense of Definition \ref{dw}.
\end{theo}
{\bf Remark.} \quad
i) \ In Keller-Segel-production systems of the form (\ref{ks}) with $D\equiv 1$ and $\chi\equiv 1$ as well as $K(u,v)=u-v$,
available existence results including
measure-type initial data seem restricted to one-dimensional and certain subcritical-mass two-dimensional settings
(\cite{biler_studiamath1995}, \cite{raczynski}, \cite{win_smoothing_1d}; cf.~also \cite{heihoff}).
Only under the influence of certain additional superlinear zero order degradation mechanisms of logistic type, 
results on instantaneous smoothing of comparably strong singularities seem to have been established in the literature
(\cite{win_AMPA}, \cite{lankeit_PROCA}).\abs
ii) \ Even for the classical chemotaxis-consumption version of (\ref{ks}) given by
\be{KSc}
	\left\{ \begin{array}{l}	
	u_t = \Del u - \na \cdot (u\na v), \\[1mm]
	v_t = \Delta v -uv,
	\end{array} \right.
\ee
%obtained on letting $D\equiv 1$ and $\chi\equiv 1$ as well as $K(u,v)=-uv$ in (\ref{ks}), 
the apparently only existence result covering large measures as initial population distributions
is limited to planar and radially symmetric frameworks (\cite{wwx_sci_chi_math}).\abs
iii) \ Upon imposing further restrictions on the regularity of $\phi$ and the initial data, higher regularity features
of the obtained solutions can be derived. Pursuing this in detail would go beyond the scope of the present study, however,
and will be addressed in \cite{genglin}.\abs
Now in the presence of slightly more regular coefficient functions $\phi$,
a second relaxation effect, unconditional with respect to the initial data from the above class and especially independent
of their size, can be observed on large time scales.
In addressing this second main objective of this study, we let $A$ denote 
the realization of $-\Del$ under homogeneous Neumann boundary conditions in 
$L^2_\perp(\Om):=\{ \vp\in L^2(\Om) \ | \ \io \vp=0\}$,
with its domain given by $W_N^{2,2}(\Om) \cap L^2_\perp(\Om)$, where for $p\in [1,\infty]$ we set 
$W_N^{2,p}(\Om):=\{\vp \in W^{2,p}(\Om) \ | \ \frac{\pa\vp}{\pa\nu}=0 \mbox{ on } \pO\}$;
furthermore, with regard to spatial averages we will refer to the notation
$\ov{\vp}:=\frac{1}{|\Om|} \vp(\one_\Om)$ for $\vp\in (C^0(\bom))^\star$, which reduces to the identity
$\ov{\vp}=\frac{1}{|\Om|} \io \vp$ whenever $\vp\in L^1(\Om)$.\abs
The second of our main results can thereby be formulated as a statement on long-term stabilization toward homogeneous states
in an appropriate topological setting:
\begin{theo}\label{theo11}
  Let $n\in\{1,2,3\}$ and $\Om\subset\R^n$ be a smoothly bounded domain, and assume that 
  \be{phi1}
	\phi\in C^1([0,\infty))
	\mbox{\quad is such that \quad } \phi(\xi)>0 \mbox{ for all } \xi\ge 0,
  \ee
  and that (\ref{init}) holds.
  Then one can find a global very weak solution $(u,v)$ of (\ref{0}) which satisfies (\ref{reg}) and for which 
  there exists a null set $N\subset (0,\infty)$ such that $A^{-\frac{1}{2}}\big(u(\cdot,t)-\ouz\big)\in L^2(\Om)$
  for all $t\in (0,\infty)\sm N$, and that
  \be{11.1}
	A^{-\frac{1}{2}}\big(u(\cdot,t)-\ouz\big)
	\to 0
	\quad \mbox{in } L^2(\Om)
	\qquad \mbox{as } (0,\infty)\sm N \ni t\to\infty
  \ee
  and
  \be{11.2}
	v(\cdot,t)\to 0
	\quad \mbox{in } L^\infty(\Om)
	\qquad \mbox{as } (0,\infty)\sm N \ni t\to\infty.
  \ee
\end{theo}
{\bf Remark.} \quad
i) \ While marking a sharp contrast to the mentioned findings on persistently non-homogeneous behavior in
corresponding chemotaxis-production versions of (\ref{ks})-(\ref{dc}) 
(\cite{fujie_jiang_CVPDE}, \cite{fujie_senba}, \cite{DLTW})), the outcome of Theorem \ref{theo11} rather parallels
known results on large time convergence to constant equilibria in the Keller-Segel consumption system
(\ref{KSc}) (\cite{taowin_consumption}). 
How far asymptotic homogenization properties of this style indeed constitute a common feature of taxis-absorption 
interaction involving non-degenerate diffusion, however, is to be addressed in forthcoming studies.\abs
{\bf Main ideas.} \quad
In line with the particular structure distinguishing (\ref{00}) from (\ref{KSc}),  
at the core of our reasoning will be a duality-based argument related to the behavior of 
\be{m1}
	\io \Big|A^{-\frac{1}{2}}\big(u-\ouz\big)\Big|^2
\ee
along trajectories in suitably regularized variants of (\ref{0}) (cf.~(\ref{0eps})).
Here in a fundamental inequality describing the evolution thereof (Lemma \ref{lem3}), in the general setting
of Theorem \ref{theo10} a fairly rough estimation of a corresponding exciting contribution is sufficient to
turn this into a quasi-energy inequality (Lemma \ref{lem4}). 
Again thanks to the favorable structure of (\ref{00}), the a priori information thereby obtained in time intervals of the
form $(\tau,\infty)$ for arbitrary $\tau>0$ can be combined with a time-independent 
$(W_N^{2,\infty}(\Om))^\star$-valued boundedness feature of both time derivatives in (\ref{0}) (Lemma \ref{lem44}),
allowing for suitably control of the solution behavior near $t=0$, to establish Theorem \ref{theo10} in Section \ref{sect3}.\abs
Section \ref{sect4} will thereafter reveal that in the low-dimensional and slightly more regular context of Theorem \ref{theo11},
the forcing contribution to the evolution of the functional in (\ref{m1}) can actually even be controlled in terms of 
suitably dissipated quantities (Lemmata \ref{lem6} and \ref{lem7}). 
An accordingly discovered energy feature will hence form the basis for our derivation of
both stabilization statements from Theorem \ref{theo11} (Lemmata \ref{lem8} and \ref{lem9}).
\mysection{Preliminaries}
To begin with, let us specify the concept of solvability to be pursued in this paper.
\begin{defi}\label{dw}
  Let $\phi\in C^0([0,\infty))$ be nonnegative, and assume that $u_0\in (C^0(\bom))^\star$ and $v_0\in L^\infty(\Om)$ are
  nonnegative.
  Then a pair $(u,v)$ of nonnegative functions
  \be{reg1}
	\left\{ \begin{array}{l}
	u\in C^0_{w-\star}([0,\infty);(C^0(\bom))^\star) \cap L^1_{loc}(\bom\times (0,\infty))
	\qquad \mbox{and} \\[1mm]
	v\in C^0_{w-\star}([0,\infty);L^\infty(\Om))
	\ear
  \ee
  will be called a {\em global very weak solution} of (\ref{0}) 
  if $u(\cdot,0)=u_0$ in $(C^0(\bom))^\star$ and $v(\cdot,0)=v_0$ in $L^\infty(\Om)$,
  and if for each $\vp\in C_0^\infty(\bom\times (0,\infty))$ fulfilling $\frac{\pa\vp}{\pa\nu}=0$ on $\pO\times (0,\infty)$ we have
  \bea{wu}
	-\int_0^\infty \io u\vp_t
	= \int_0^\infty \io u\phi(v) \Del\vp
  \eea
  and
  \be{wv}
	- \int_0^\infty \io v\vp_t  
	= \int_0^\infty \io v\Del\vp
	- \int_0^\infty \io uv\vp.
  \ee
\end{defi}
In order to construct such solutions to (\ref{0}) as limits of solutions to suitably regularized problems,
we approximate the motility function $\phi$ as well as initial data $(u_0,v_0)$ in (\ref{0}) by introducing 
families of functions $(\pe)_{\eps\in (0,1)}$, $(u_{0\eps})_{\eps\in (0,1)}$ and $(v_{0\eps})_{\eps\in (0,1)}$ with the properties that 
%
%
	%$(\phi_\eps)_{\eps\in (0,1)} \subset C^3([0,\infty))$ be such that 
\be{pe}
	\lbal
	(\pe)_{\eps\in (0,1)} \subset C^3([0,\infty)) \mbox{\quad is such that } \\[1mm] 
	\pe\ge \phi \mbox{ on $[0,\infty)$ for all } \eps\in (0,1), \quad \mbox{and that} \\[1mm]
	\pe\to\phi \mbox{ in } C^0_{loc}([0,\infty))
	\quad \mbox{as } \eps\searrow 0,
	\ear
\ee
and that moreover
\be{ie}
	\lbal
	(u_{0\eps})_{\eps\in (0,1)} \subset W^{1,\infty}(\Om) \mbox{ and } (v_{0\eps})_{\eps\in (0,1)} \subset W^{1,\infty}(\Om)
		\quad \mbox{are such that} \\[1mm]
	u_{0\eps}\ge 0 \ \mbox{ and } \ v_{0\eps} > 0 \quad \mbox{ in $\bom$ \qquad for all $\eps\in (0,1)$, \qquad that} \\[1mm]
	\io u_{0\eps}= \ouz |\Om|
		\quad \mbox{and} \ \|v_{0\eps}\|_{L^\infty(\Om)} \le \|v_0\|_{L^\infty(\Om)}+1
	\quad \mbox{for all $\eps\in (0,1)$, \quad and that} \\[1mm]
	u_{0\eps} \wsto u_0 \mbox{ in } (C^0(\bom))^\star \quad \mbox{and} \quad
	v_{0\eps} \wsto v_0 \mbox{ in } L^\infty(\Om)
	\quad \mbox{as } \eps\searrow 0.
	\ear
\ee

For $\eps\in (0,1)$, we then consider the regularized variant of (\ref{0}) given by
\be{0eps}
	\lball
	u_{\eps t} = \Del \big(\ueps\phi_\eps(\veps)\big),
	\qquad & x\in\Omega, \ t>0, \\[1mm]
	v_{\eps t} = \Delta\veps - \frac{\ueps\veps}{1+\eps\ueps},
	\qquad & x\in\Omega, \ t>0, \\[1mm]
	\frac{\partial\ueps}{\partial\nu}=\frac{\partial\veps}{\partial\nu}=0,
	\qquad & x\in\pO, \ t>0, \\[1mm]
	\ueps(x,0)=u_{0\eps}(x), \quad
	\veps(x,0)=v_{0\eps}(x),
	\qquad & x\in\Omega，
	\ear
\ee
which is globally solved in the classical sense:
\begin{lem}\label{lem_loc}
  For each $\eps\in (0,1)$ there exist 
  \be{01.1}
	\lbal
	\ueps\in C^0(\bom\times [0,\infty)) \cap C^{2,1}(\bom\times (0,\infty))
	\qquad \mbox{and} \\[1mm]
	\veps \in \bigcap_{q>2} C^0([0,\infty);W^{1,q}(\Om)) \cap C^{2,1}(\bom\times (0,\infty))
	\ear
  \ee
  such that $\ueps \ge 0$ and $\veps>0$ in $\bom\times [0,\infty)$, 
  and that $(\ueps,\veps)$ solves (\ref{0eps}) in the classical sense. Furthermore, the solution satisfies 
  \be{mass}
  	\io \ueps(\cdot,t) = \io u_0
  	\qquad \mbox{for all $t>0$ and } \eps\in (0,1),
    \ee
    and 
    \be{vinfty}
  	\|\veps(\cdot,t)\|_{L^\infty(\Om)} \le \|\veps(\cdot,t_0)\|_{L^\infty(\Om)}
  	\qquad \mbox{for all $t_0\ge 0$, $t>t_0$ and } \eps\in (0,1).
  	\ee
\end{lem}
\proof
  We start by asserting the local classical solvability for (\ref{0eps}) by means of the well-established
  parabolic theory from (\cite{amann}). To this end, we fix $\delta_0>0$, and for $\eps\in(0,1)$ introducing
  $D_0:= (0, \infty)\times(-\delta_0, \infty)$ as well as
  \begin{equation*}
  	A_{\eps}\Big( \begin{array}{c} \eta \\ \xi \end{array} \Big)
  	:=\left( \begin{array}{cc}
  	               1 & 0\\
  	\xi\phi_\eps'(\eta)   & \phi_\eps(\eta)
  	\end{array} \right)
  	\quad \mbox{and} \qquad
	f_{\eps}\Big( \begin{array}{c} \eta \\ \xi \end{array} \Big)
  	:=\left( \begin{array}{c}
  	-\frac{\xi\eta}{1+\eps\xi}
	\\ 
	0
  	\end{array} \right)
  	\quad \mbox{for } \ \Big( \begin{array}{c} \eta \\ \xi \end{array} \Big) \in D_0,
  \end{equation*}
  we may recast (\ref{0eps}) as the quasilinear problem
  \be{00eps}
  \left\{ \begin{array}{lcll}
  	Z_{\eps t} &=& \na\cdot\big( A_{\eps}(Z_\eps) \na Z_{\eps }\big)
  	+  f_{\eps}(Z_{\eps}) ,
  	& x\in\Omega, \ t>0, \\[2mm]
  	& & \hspace*{-15mm}
  	\frac{\pa Z_{\eps }}{\pa  \nu}=0,
  	& x\in\pO, \ t>0, \\[1mm]
  	& & \hspace*{-15mm}
  	Z_\eps(x,0)=\left( \begin{array}{c} v_{0\eps}(x) \\ u_{0\eps}(x) \end{array} \right),
  	& x\in\Omega,
  \end{array} \right.
  \ee 
  where $Z_{\eps}=(\veps,\ueps)$. Using (\ref{phi}) and (\ref{pe}), we observe that for each $U\in D_0$,
  $A_{\eps}(U)$ is a positive definite matrix of lower triangular form, so that from (\cite[Theorem 1]{amann}),
  in line with (\ref{ie}) we deduce that there exists $\tme\in(0,\infty]$ such that (\ref{0eps}) possesses a classical solution
  $(\ueps,\veps)$ which is such that $\veps>0$ in $\bom\times[0,\tme)$, and that
  \be{tme}
 	\mbox{if $\tme<\infty$}, \quad then \quad 
  	\limsup_{t\nearrow \tme} \|\ueps(\cdot,t)\|_{ L^{\infty}(\Omega)} =\infty.
  \ee
  Moreover, taking into account $u_{0\eps}\ge 0$ in (\ref{ie}), we obtain $\ueps\ge 0$ in $\bom\times[0,\tme)$ through
  a simple comparison argument. Next after an integration performed in the first equation in (\ref{0eps}), we find that
  \be{mass1}
  	\io \ueps(\cdot,t) = \io u_0
  	\qquad \mbox{for all $t\in (0,\tme)$ and } \eps\in (0,1),
  \ee
  whereas an application of the maximum principle to the second equation therein shows that
  \be{vinfty1}
  	\|\veps(\cdot,t)\|_{L^\infty(\Om)} \le \|\veps(\cdot,t_0)\|_{L^\infty(\Om)}
  	\qquad \mbox{for all $t_0\in [0,\tme)$, $t\in (t_0,\tme)$ and } \eps\in (0,1).
  \ee
  To finally prove that $\tme=\infty$, assuming on the contrary that $\tme<\infty$ for some $\eps\in(0,1)$, we use
  (\ref{vinfty1}) to see that $\Big\|\frac{\ueps(\cdot,t)\veps(\cdot,t)}{1+\eps\ueps(\cdot,t)}\Big\|_{L^{\infty}(\Om)}
  \le\frac{\|v_{0\eps}\|_{L^{\infty}(\Om)}}{\eps}$ for all $t\in (0,\tme)$, and can then rely on 
  standard parabolic regularity theory 		%(\cite{horstmann_winkler}) 
  applied to 
  the second equation in (\ref{0eps}) to find $c_1(\eps)>0$ such that 
  \bea{l1}
  	\|\na \veps(\cdot,t)\|_{L^{\infty}(\Om)} \le c_1(\eps)
  	\qquad \mbox{for all $t\in (0,\tme)$}.
  \eea
  In addition, again thanks to (\ref{vinfty1}), (\ref{phi}) and (\ref{pe}) warrant the existence of positive constants $c_2$ 
  and $c_3(\eps)$ fulfilling
  \bea{l2}
 	\phi_\eps(\veps)\ge c_2\quad \mbox{and}\quad
 	\frac{(\phi_\eps'(\veps))^2}{\phi_\eps(\veps)}\le c_3(\eps)
	\qquad \mbox{in } \Om\times (0,\tme).
  \eea 
  Therefore, integrating by parts in the first equation from (\ref{0eps}) and using the Cauchy-Schwarz inequality and
  (\ref{l1}), we see that whenever $p>1$,
  \bea{l3}
 	\frac{d}{dt}\io \ueps^p 
 	&=& -p(p-1)\io \phi_\eps(\veps)\ueps^{p-2}|\na \ueps|^2
          - p(p-1)\io \phi_\eps'(\veps)\ueps^{p-1}\na \ueps\cdot\na\veps\nn\\
 	&\le& -\frac{p(p-1)}{2}\io \phi_\eps(\veps)\ueps^{p-2}|\na \ueps|^2
         + \frac{p(p-1)}{2}\io \frac{(\phi_\eps'(\veps))^2}{\phi_\eps(\veps)}\ueps^p|\na \veps|^2\nn\\
 	&\le& -\frac{2c_2(p-1)}{p}\io|\na \ueps^{\frac{p}{2}}|^2
         + \frac{c_1^2(\eps)c_3(\eps)p(p-1)}{2}\io \ueps^p
  	\qquad \mbox{for all $t\in (0,\tme)$}.      
  \eea
  According to an Ehrling type inequality associated with the compactness of the embedding $W^{1,2}(\Om)\hra L^2(\Om)$,
  from (\ref{mass1}) it follows that for any such $p$, with some $c_4(p,\eps)>0$ we have
  \bas
  	\frac{c_1^2(\eps)c_3(\eps)p(p-1)}{2}\io \ueps^p
   	&=& \frac{c_1^2(\eps)c_3(\eps)p(p-1)}{2} \|\ueps^{\frac{p}{2}}\|_{L^2(\Om)}^2\\
  	&\le& \frac{2c_2(p-1)}{p} \|\na\ueps^{\frac{p}{2}}\|_{L^2(\Om)}^2
         + c_4(p,\eps) \|\ueps^{\frac{p}{2}}\|_{L^{\frac{2}{p}}(\Om)}^2\\
   	&=&  \frac{2c_2(p-1)}{p} \|\na\ueps^{\frac{p}{2}}\|_{L^2(\Om)}^2
         + c_4(p,\eps) \Big\{\io u_0\Big\}^p 
  	\qquad \mbox{for all $t\in (0,\tme)$},
  \eas
  which, when substituted back into (\ref{l3}), upon an integration in time entails that
  \bas
  	\io \ueps^p 
  	&\le& c_4(p,\eps) \Big\{\io u_0\Big\}^p \tme + \io u_{0\eps}^p\\
  	&\le& c_4(p,\eps) \Big\{\io u_0\Big\}^p \tme + \|u_{0\eps}\|_{W^{1,\infty}(\Om)}^p|\Om|
  	\qquad \mbox{for all $t\in (0,\tme)$}.
  \eas
  Since $p>1$ was arbitrary here, in view of a standard Moser-type iteration (\cite{taowin_subcrit}) 
  we may find $c_5(\eps)>0$ such that
  \bas
 	\|\ueps(\cdot,t)\|_{L^{\infty}(\Omega)} \le c_5(\eps)
  	\qquad \mbox{for all $t\in (0,\tme)$},
  \eas
  which contradicts (\ref{tme}) and thus confirms that, indeed, $\tme=\infty$. 
\qed
\mysection{Instantaneous relaxation. Proof of Theorem \ref{theo10}}\label{sect3}
The following basic properties of solutions to (\ref{0eps}) can be established in a straightforward manner.
\begin{lem}\label{lem1}
  Let $n\ge 1$, and assume (\ref{ie}) and (\ref{pe}).
   Then 
   \begin{comment}
  \be{mass}
	\io \ueps(\cdot,t) = \io u_0
	\qquad \mbox{for all $t>0$ and } \eps\in (0,1),
  \ee
  and we have
  \be{vinfty}
	\|\veps(\cdot,t)\|_{L^\infty(\Om)} \le \|\veps(\cdot,t_0)\|_{L^\infty(\Om)}
	\qquad \mbox{for all $t_0\ge 0$, $t>t_0$ and } \eps\in (0,1).
  \ee
  Moreover,
  \end{comment}
  \be{uv}
	\int_0^\infty \io \frac{\ueps\veps}{1+\eps\ueps} \le |\Om| \cdot \big(\|v_0\|_{L^\infty(\Om)}+1\big)
	\qquad \mbox{for all } \eps\in (0,1)
  \ee
  and
  \be{gradv}
	\int_0^\infty \io |\na\veps|^2 \le \frac{1}{2} |\Om| \cdot \big( \|v_0\|_{L^\infty(\Om)}+1\big)^2
	\qquad \mbox{for all } \eps\in (0,1).
  \ee
\end{lem}
\proof
  For $p\ge 1$, we integrate using the second equation in (\ref{0eps}) to see that
  \be{1.1}
	\frac{1}{p} \frac{d}{dt} \io \veps^p + (p-1) \io \veps^{p-2} |\na\veps|^2
	+ \io \frac{\ueps \veps^p}{1+\eps\ueps} =0
	\qquad \mbox{for all $t>0$ and } \eps\in (0,1).
  \ee
  When specified to the case $p=1$, upon a time integration this implies that due to (\ref{ie}),
  \bas
	\int_0^T \io \frac{\ueps\veps}{1+\eps\ueps} 
	= - \io \veps(\cdot,T)+\io v_{0\eps}
	\le |\Om| \cdot \big(\|v_0\|_{L^\infty(\Om)}+1\big)
	\qquad \mbox{for all $T>0$ and } \eps\in (0,1),
  \eas
  and that thus (\ref{uv}) holds.
  Secondly, in the case when $p=2$, an integration of (\ref{1.1}) shows that 
  for all $T>0$ and $\eps\in (0,1)$,
  \bas
	\int_0^T \io |\na\veps|^2
	&=& - \int_0^T \io \frac{\ueps\veps^2}{1+\eps\ueps}
	- \frac{1}{2} \io \veps^2(\cdot,T)
	+ \frac{1}{2} \io v_{0\eps}^2
	\le \frac{1}{2} |\Om| \big(\|v_0\|_{L^\infty(\Om)}+1\big)^2,
  \eas
  and thereby establishes (\ref{gradv}).
\qed
\begin{comment}
  Finally, (\ref{mass}) immediately results from an integration in (\ref{0eps}) again using (\ref{ie}), while
  (\ref{vinfty}) is a consequence of the maximum principle.
\end{comment}

%
%
%
%
The following outcome of an essentially straightforward testing procedure is formulated in such a way that
it can not only serve as an ingredient in our derivation of boundedness features in the general setting of Theorem \ref{theo10},
but also be used in the course of our asymptotic analysis in the low-dimensional situations addressed by Theorem \ref{theo11}.
\begin{lem}\label{lem2}
  If $n\ge 1$, and if (\ref{ie}) and (\ref{pe}) hold, then one can find $\Gamma_1>0$ such that
  \be{2.1}
	\frac{d}{dt} \io |\na\veps|^2 
	+ \frac{1}{2} \io |\Del\veps|^2
	+ \frac{1}{\Gamma_1} \io |\na\veps|^4
	\le \Gamma_1 \io (\ueps-\ouz)^2
	\qquad \mbox{for all $t>0$ and } \eps\in (0,1).
  \ee
\end{lem}
\proof
  Let $F_\eps(\xi):=\frac{\xi}{1+\eps\xi}$ for $\xi\ge 0$ and $\eps\in (0,1)$. Then since $0\le F_\eps'\le 1$ for all
  $\eps\in (0,1)$, from the mean value theorem it follows that 
  \bas
	\Big|\frac{\ueps}{1+\eps\ueps} - \frac{\ouz}{1+\eps\ouz}\Big|
	= \big|F_\eps(\ueps)-F_\eps(\ouz)\big|
	\le |\ueps-\ouz|
	\quad \mbox{in } \Om\times (0,\infty)
	\qquad \mbox{for all } \eps\in (0,1).
  \eas
  Therefore, if we multiply the second equation in (\ref{0eps}) by $-\Del\veps$ and integrate by parts and using Young's inequality,
  (\ref{vinfty}) and (\ref{ie}), we see that
  \bea{2.2}
	\frac{1}{2} \frac{d}{dt} \io |\na\veps|^2
	+ \io |\Del\veps|^2
	&=& \io \frac{\ueps\veps}{1+\eps\ueps} \Del\veps \nn\\
	&\le& \frac{1}{2} \io |\Del\veps|^2 
	+ \frac{1}{2} \io \Big( \frac{\ueps}{1+\eps\ueps} - \frac{\ouz}{1+\eps\ouz}\Big)^2 \veps^2 \nn\\
	&\le& \frac{1}{2} \io |\Del\veps|^2 
	+ \frac{1}{2} \cdot \big( \|v_0\|_{L^\infty(\Om)}+1\big)^2 \io (\ueps-\ouz)^2
%	\qquad \mbox{for all $t>0$ and } \eps\in (0,1),
  \eea
  for all $t>0$ and $\eps\in (0,1)$,
  because
  \bas
	\io \veps\Del\veps = - \io |\na\veps|^2 \le 0
	\qquad \mbox{for all $t>0$ and } \eps\in (0,1).
  \eas
  Now from elliptic regularity theory we obtain that with some $c_1>0$ we have
  \bas
	\io |D^2\vp|^2 \le c_1 \io |\Del\vp|^2
	\qquad \mbox{for all $\vp\in C^2(\bom)$ fulfilling $\frac{\pa\vp}{\pa\nu}=0$ on } \pO,
  \eas
  so that another integration by parts together with the Cauchy-Schwarz inequality shows that
  \bas
	\io |\na\veps|^4
	&=& - \io \veps \cdot \Big\{ 2\na\veps\cdot (D^2\veps\cdot\na\veps) + |\na\veps|^2 \Del\veps \Big\} \\
	&\le& 2 \io \veps |\na\veps|^2 |D^2\veps|
	+ \io \veps |\na\veps|^2 |\Del\veps| \\
	&\le& \|\veps\|_{L^\infty(\Om)} \cdot \bigg\{ \io |\na\veps|^4 \bigg\}^\frac{1}{2} \cdot \Bigg\{
	2\cdot\bigg\{ \io |D^2\veps|^2 \bigg\}^\frac{1}{2} 
	+ \bigg\{ \io |\Del\veps|^2 \bigg\}^\frac{1}{2} \Bigg\} \\
	&\le& \|\veps\|_{L^\infty(\Om)} \cdot \bigg\{ \io |\na\veps|^4 \bigg\}^\frac{1}{2}
	\cdot (2\sqrt{c_1} +1) \cdot \bigg\{ \io |\Del\veps|^2 \bigg\}^\frac{1}{2}
	\qquad \mbox{for all $t>0$ and } \eps\in (0,1),
  \eas
  and thus, by (\ref{vinfty}) and (\ref{ie}),
  \bas
	\io |\na\veps|^4
	\le \|\veps\|_{L^\infty(\Om)}^2 (2\sqrt{c_1}+1)^2 \io |\Del\veps|^2
	\le c_2 \io |\Del\veps|^2
	\qquad \mbox{for all $t>0$ and } \eps\in (0,1)
  \eas
  with $c_2:=(2\sqrt{c_1}+1)^2 (\|v_0\|_{L^\infty(\Om)}+1)^2$.
  Therefore, (\ref{2.2}) entails that
  \bas
	\frac{1}{2} \frac{d}{dt} \io |\na\veps|^2
	+ \frac{1}{4} \io |\Del\veps|^2
	+ \frac{1}{4c_2} \io |\na\veps|^4
	&\le& \frac{1}{2} \frac{d}{dt} \io |\na\veps|^2 + \frac{1}{2} \io |\Del\veps|^2 \\
	&\le& \frac{1}{2} \big(\|v_0\|_{L^\infty(\Om)}+1\big)^2 \io (\ueps-\ouz)^2
%	\qquad \mbox{for all $t>0$ and } \eps\in (0,1),
  \eas
  for all $t>0$ and $\eps\in (0,1)$, and that hence (\ref{2.1}) holds with
  $\Gamma_1:=\max\{2c_2 \, , \, (\|v_0\|_{L^\infty(\Om)}+1)^2\}$.
\qed
Also our first information on an evolution property of the first solution components,
acting in an $H^{-1}$ framework reminiscent of that already resorted to in \cite{taowin_M3AS},
is at this stage kept general enough so as to remain 
compatible with our analysis of both Theorem \ref{theo10} and Theorem \ref{theo11}.
\begin{lem}\label{lem3}
  Suppose that $n\ge 1$, and that (\ref{ie}) and (\ref{pe}) are satisfied. Then there exists $\Gamma_2>0$ with the property that
  \be{3.1}
	\frac{d}{dt} \io \Big| A^{-\frac{1}{2}}(\ueps-\ouz)\Big|^2
	+ \frac{1}{\Gamma_2} \io (\ueps-\ouz)^2
	\le \Gamma_2 \io \Big| \ouz \pe(\veps)- \ov{\ueps(\cdot,t)\pe\big(\veps(\cdot,t)\big)} \Big|^2
%	\qquad \mbox{for all $t>0$ and } \eps\in (0,1).
  \ee
  for all $t>0$ and $\eps\in (0,1)$.
\end{lem}
\proof
  According to (\ref{0eps}), we have
  \bas
	\pa_t (\ueps-\ouz) = -A \big( \ueps\pe(\veps)-\ov{\ueps\pe(\veps)}\big)
	\quad \mbox{in } \Om\times (0,\infty)
	\qquad \mbox{for all } \eps\in (0,1),
  \eas
  which when tested against $A^{-1} (\ueps-\ouz)$ implies that since both $A^{-\frac{1}{2}}$ and $A^{-1}$ are self-adjoint,
  \bas
	\frac{1}{2} \frac{d}{dt} \io \Big|A^{-\frac{1}{2}}(\ueps-\ouz)\Big|^2
	&=& - \io \big(\ueps\pe(\veps) - \ov{\ueps\pe(\veps)}\big) \cdot (\ueps-\ouz) \\
	&=& - \io \big(\ueps\pe(\veps)-\ouz\pe(\veps) + \ouz \pe(\veps) - \ov{\ueps\pe(\veps)} \big) \cdot (\ueps-\ouz) \\
	&=& - \io (\ueps-\ouz)^2 \pe(\veps)
	- \io \big(\ouz \pe(\veps) - \ov{\ueps\pe(\veps)}\big)\cdot (\ueps-\ouz)
  \eas
  for all $t>0$ and $\eps\in (0,1)$.
  Since (\ref{pe}) together with (\ref{phi}), (\ref{vinfty}) and (\ref{ie}) entails the existence of $c_1>0$ such that
  $\pe(\veps)\ge c_1$ in $\Om\times (0,\infty)$ for all $\eps\in (0,1)$, and since
  \bas
	- \io \big(\ouz \pe(\veps) - \ov{\ueps\pe(\veps)}\big)\cdot (\ueps-\ouz)
	\le \frac{c_1}{2} \io (\ueps-\ouz)^2
	+ \frac{1}{2c_1} \io \Big| \ouz\pe(\veps)-\ov{\ueps\pe(\veps)}\Big|^2
%	\qquad \mbox{for all $t>0$ and } \eps\in (0,1)
  \eas
  for all $t>0$ and $\eps\in (0,1)$
  thanks to Young's inequality, this already leads to (\ref{3.1}) if we let $\Gamma_2:=\frac{1}{c_1}$, for instance.
\qed
In our first application of Lemma \ref{lem2} and Lemma \ref{lem3}, we may estimate the expression on the right
of (\ref{3.1}) in a fairly rough manner.
Actually focusing especially on initial relaxation features here, we can thereby, after all, make sure that
a functional of the form in (\ref{m1}) satisfies an ODI containing a superlinear absorptive term, an appropriate exploitation 
of which yields a priori information within time intervals of the form $(\tau,\infty)$ for arbitrary $\tau>0$:
\begin{lem}\label{lem4}
  Let $n\ge 1$, and assume (\ref{ie}) and (\ref{pe}).
  Then for all $\tau>0$ there exists $C(\tau)>0$ such that
  \be{4.1}
	\io \Big| A^{-\frac{1}{2}}\big(\ueps(\cdot,t)-\ouz\big)\Big|^2
	\le C(\tau)
	\qquad \mbox{for all $t>\tau$ and } \eps\in (0,1),
  \ee
  that
  \be{4.01}
	\io \big|\na\veps(\cdot,t)\big|^2 
	\le C(\tau)
	\qquad \mbox{for all $t>\tau$ and } \eps\in (0,1),
  \ee
  that
  \be{4.2}
	\int_t^{t+1} \io \ueps^2
	\le C(\tau)
	\qquad \mbox{for all $t>\tau$ and } \eps\in (0,1),
  \ee
  that
  \be{4.3}
	\int_t^{t+1} \io |\Del\veps|^2
	\le C(\tau)
	\qquad \mbox{for all $t>\tau$ and } \eps\in (0,1),
  \ee
  and that
  \be{4.4}
	\int_t^{t+1} \io |\na\veps|^4
	\le C(\tau)
	\qquad \mbox{for all $t>\tau$ and } \eps\in (0,1).
  \ee
\end{lem}
\proof
  As a consequence of (\ref{pe}), (\ref{phi}), (\ref{vinfty}) and (\ref{ie}), we can find $c_1>0$ such that $\pe(\veps)\le c_1$ in
  $\Om\times (0,\infty)$ for all $\eps\in (0,1)$. Since thus especially
  \bas
	\ov{\ueps\pe(\veps)} 
	= \frac{1}{|\Om|} \io \ueps\pe(\veps) 
	\le \frac{c_1}{|\Om|} \io \ueps= c_1 \ouz
	\qquad \mbox{for all $t>0$ and } \eps\in (0,1)
  \eas
  by (\ref{mass}), on the right-hand side of (\ref{3.1}) we can estimate
  \bea{4.5}
	\Gamma_2 \io \Big| \ouz \pe(\veps)- \ov{\ueps(\cdot,t)\pe\big(\veps(\cdot,t)\big)} \Big|^2
	&\le& \Gamma_2 \io |c_1\ouz+c_1\ouz|^2 \nn\\
	&=& c_2:=4c_1^2 \Gamma_2 \ouz^2 |\Om|
	\qquad \mbox{for all $t>0$ and } \eps\in (0,1).
  \eea
  From Lemma \ref{lem3} and Lemma \ref{lem2} we therefore obtain that if we let $a:=\frac{1}{4\Gamma_1 \Gamma_2}$,
  with $\Gamma_1$ taken from Lemma \ref{lem2}, and define
  \be{4.6}
	\yeps(t):=\io \Big|A^{-\frac{1}{2}}\big(\ueps(\cdot,t)-\ouz\big)\Big|^2
	+ a \io \big|\na\veps(\cdot,t)\big|^2,
	\qquad t>0, \ \eps\in (0,1),
  \ee
  as well as
  \be{4.7}
	\geps(t):=\frac{1}{2\Gamma_2} \io \big(\ueps(\cdot,t)-\ouz\big)^2
	+ \frac{a}{2} \io \big|\Del\veps(\cdot,t)\big|^2
 	+ \frac{a}{2\Gamma_1} \io \big|\na\veps(\cdot,t)\big|^4,
	\qquad t>0, \ \eps\in (0,1),
  \ee
  then
  \bea{4.8}
	& & \hs{-20mm}
	\yeps'(t) + \geps(t)
	+ \frac{1}{4\Gamma_2} \io (\ueps-\ouz)^2
	+ \frac{a}{2\Gamma_1} \io |\na\veps|^4 \nn\\
	&\le& \bigg\{ -\frac{1}{4\Gamma_2} \io (\ueps-\ouz)^2 + c_2 \bigg\} \nn\\
	& & + a\cdot \Gamma_1 \io (\ueps-\ouz)^2  \nn\\[2mm]
	&=& c_2
	\qquad \mbox{for all $t>0$ and } \eps\in (0,1).
  \eea
  In order to turn this into a superlinearly damped ODI for $(\yeps)_{\eps\in (0,1)}$, we pick any $\beta>\frac{1}{2}$ such that
  $\beta>\frac{n}{4}$ and note that then the inequalities $0<\frac{1}{2}<\beta$ enable us to invoke a standard interpolation
  inequality (\cite[Theorem 14.1]{friedman} to find $c_3>0$ fulfilling
  \bas
	\|A^{-\frac{1}{2}} \vp\|_{L^2(\Om)}
	\le c_3 \|\vp\|_{L^2(\Om)}^\theta \|A^{-\beta}\vp\|_{L^2(\Om)}^{1-\theta}
	\qquad \mbox{for all } \vp\in L^2_\perp(\Om)
  \eas
  with $\theta:=\frac{2\beta-1}{2\beta}\in (0,1)$.
  Since the restriction $\beta>\frac{n}{4}$ warrants that for the corresponding domains of definition we have 
  $D(A^\beta) \hra L^\infty(\Om)$ and hence $L^1(\Om) \hra D(A^{-\beta})$ (\cite{henry}), we may combine this with (\ref{mass})
  to see that with some $c_4>0$ we have
  \bas
	\bigg\{ \io \Big|A^{-\frac{1}{2}}(\ueps-\ouz)\Big|^2\bigg\}^\frac{1}{\theta} \le c_4 \io (\ueps-\ouz)^2
	\qquad \mbox{for all $t>0$ and } \eps\in (0,1).
  \eas
  Writing $\kappa:=\min\{ \frac{1}{\theta} \, , \, 2 \}>1$ and
  \be{4.10}
	c_5:=\min \bigg\{ \frac{1}{2^{\kappa+1} c_4 \Gamma_2} \, , \, \frac{1}{2^\kappa \Gamma_1 a^{\kappa-1} |\Om|} \bigg\},
  \ee
  we thus infer that in line with (\ref{4.6}) and thanks to the Cauchy-Schwarz inequality,
  \bas
	c_5 \yeps^\kappa(t)
	&\le& 2^{\kappa-1} c_5 \cdot \bigg\{ \io \Big|A^{-\frac{1}{2}}(\ueps-\ouz)\Big|^2 \bigg\}^\kappa
	+ 2^{\kappa-1} a^\kappa c_5 \cdot \bigg\{ \io |\na\veps|^2 \bigg\}^\kappa \\
	&\le& 2^{\kappa-1} c_5 \cdot \Bigg\{ \bigg\{ \io \Big|A^{-\frac{1}{2}}(\ueps-\ouz)\Big|^2 \bigg\}^\frac{1}{\theta} 
		+1 \Bigg\}
	+ 2^{\kappa-1} a^\kappa c_5 \cdot \Bigg\{ \bigg\{ \io |\na\veps|^2 \bigg\}^2 +1 \Bigg\} \\
	&\le& 2^{\kappa-1} c_5 \cdot \bigg\{ c_4 \io (\ueps-\ouz)^2 + 1 \bigg\}
	+ 2^{\kappa-1} a^\kappa c_5 \cdot \bigg\{ |\Om| \io |\na\veps|^4 + 1 \bigg\} \\
	&\le& \frac{1}{4\Gamma_2} \io (\ueps-\ouz)^2
	+ \frac{a}{2\Gamma_1} \io |\na\veps|^4
	+ 2^{\kappa-1} (1+a^\kappa) c_5
	\qquad \mbox{for all $t>0$ and } \eps\in (0,1),
  \eas
  because $2^{\kappa-1} c_5 c_4 \le \frac{1}{4\Gamma_2}$ and $2^{\kappa-1} a^\kappa c_5 |\Om| \le \frac{a}{2\Gamma_1}$
  due to (\ref{4.10}).
  Consequently, (\ref{4.8}) implies that with $c_6:=c_2+2^{\kappa-1} (1+a^\kappa)c_5$ we have
  \be{4.11}
	\yeps'(t) + c_5 \yeps^\kappa(t) + \geps(t) \le c_6
	\qquad \mbox{for all $t>0$ and } \eps\in (0,1),
  \ee
  so that since for fixed $\tau>0$, 
  \bas
	\oy(t):=c_7\cdot \Big(t-\frac{\tau}{2}\Big)^{-\frac{1}{\kappa-1}} + c_7,
	\qquad t>\frac{\tau}{2},
  \eas
  with  
  \bas	
	c_7:=\max \bigg\{ \big( (\kappa-1)c_5\big)^{-\frac{1}{\kappa-1}} \, , \, \Big(\frac{c_6}{c_5}\Big)^\frac{1}{\kappa} \bigg\},
  \eas
  satisfies $\oy(t)\nearrow +\infty$ as $t\searrow\frac{\tau}{2}$ and
  \bas
	\oy'(t) + c_5 \oy^\kappa(t) + \geps(t) - c_6
	&\ge& \oy'(t) + c_5 \oy^\kappa(t) - c_6 \\
	&=& - \frac{c_7}{\kappa-1} \Big(t-\frac{\tau}{2}\Big)^{-\frac{1}{\kappa-1}-1}
	+ c_5 c_7^\kappa \cdot \bigg\{ \Big(t-\frac{\tau}{2}\Big)^{-\frac{1}{\kappa-1}} +1 \bigg\}^\kappa -c_6 \\
	&\ge& - \frac{c_7}{\kappa-1} \Big(t-\frac{\tau}{2}\Big)^{-\frac{1}{\kappa-1}-1}
	+ c_5 c_7^\kappa \Big(t-\frac{\tau}{2}\Big)^{-\frac{\kappa}{\kappa-1}} 
	+ c_5 c_7^\kappa -c_6 \\
	&=& c_5 c_7 \cdot \Big\{ c_7^{\kappa-1} - \frac{1}{(\kappa-1)c_5}\Big\} 
		\cdot \Big(t-\frac{\tau}{2}\Big)^{-\frac{\kappa}{\kappa-1}}
	+ c_5 \cdot \Big\{ c_7^\kappa-\frac{c_6}{c_5}\Big\} \\[2mm]
	&\ge& 0
	\qquad \mbox{for all } t>\frac{\tau}{2},
  \eas
  an ODE comparison argument applied to (\ref{4.11}) shows that $\yeps(t)\le \oy(t)$ for all $t>\frac{\tau}{2}$ and $\eps\in (0,1)$,
  and that hence, in particular, 
  \be{4.12}
	\yeps(t) \le c_8(\tau):=c_7\cdot \Big(\frac{\tau}{2}\Big)^{-\frac{1}{\kappa-1}} +c_7
	\qquad \mbox{for all $t>\tau$ and } \eps\in (0,1).
  \ee
  Thereupon, by direct integration in (\ref{4.11}) we obtain that
  \be{4.13}
	\int_t^{t+1} \geps(s) ds
	\le \yeps(t) + \int_t^{t+1} c_6 ds
	\le c_8(\tau) + c_6
	\qquad \mbox{for all $t>\tau$ and } \eps\in (0,1),
  \ee
  so that in view of (\ref{4.6}) and (\ref{4.7}) we infer (\ref{4.1}) and (\ref{4.01}) from (\ref{4.12}) and
  (\ref{4.2})-(\ref{4.4}) from (\ref{4.13}) if we choose $C(\tau)$ appropriately large.
\qed
A straightforward estimation of corresponding time derivatives does not only pave the way toward an Aubin-Lions type 
compactness argument, but beyond this also prepares our analysis of the solution behavior near the initial instant.
\begin{lem}\label{lem44}
  Suppose that $n\ge 1$, and that (\ref{ie}) and (\ref{pe}) hold. 
  Then one can find $C>0$ such that
  \be{44.1}
	\big\| u_{\eps t}(\cdot,t)\big\|_{(W_N^{2,\infty}(\Om))^\star} \le C
	\qquad \mbox{for all $t>0$ and } \eps\in (0,1)
  \ee
  and
  \be{44.2}
	\big\| v_{\eps t}(\cdot,t)\big\|_{(W_N^{2,\infty}(\Om))^\star} \le C
	\qquad \mbox{for all $t>0$ and } \eps\in (0,1).
  \ee
\end{lem}
\proof
  We fix $\vp\in W_N^{2,\infty}(\Om)$ and use (\ref{0eps}) to see that
  \bea{44.3}
	\bigg| \io u_{\eps t} \vp \bigg|
	&=& \bigg| \io \Del\big(\ueps\pe(\veps)\big) \vp \bigg| \nn\\
	&=& \bigg| \io \ueps\pe(\veps) \Del\vp \bigg| \nn\\[2mm]
	&\le& \|\ueps\|_{L^1(\Om)} \big\|\pe(\veps)\big\|_{L^\infty(\Om)} \|\Del\vp\|_{L^\infty(\Om)}
	\qquad \mbox{for all $t>0$ and } \eps\in (0,1)
  \eea
  as well as
  \bea{44.4}
	\bigg| \io v_{\eps t} \vp \bigg|
	&=& \bigg| \io \Del\veps \vp - \io \frac{\ueps\veps}{1+\eps\ueps} \vp \bigg| \nn\\
	&=& \bigg| \io \veps\Del\vp - \io \frac{\ueps\veps}{1+\eps\ueps} \vp \bigg| \nn\\[2mm]
	&\le& \|\veps\|_{L^1(\Om)} \|\Del\vp\|_{L^\infty(\Om)} 
	+ \|\ueps\|_{L^1(\Om)} \|\veps\|_{L^\infty(\Om)} \|\vp\|_{L^\infty(\Om)},
  \eea
  because $0\le\frac{\ueps}{1+\eps\ueps} \le \ueps$ in $\Om\times (0,\infty)$ for all $\eps\in (0,1)$.
  Since (\ref{mass}), (\ref{vinfty}), (\ref{ie}), (\ref{pe}) and (\ref{phi}) guarantee boundedness of
  $(\ueps)_{\eps\in (0,1)}$ and $(\veps)_{\eps\in (0,1)}$ in $L^\infty((0,\infty);L^1(\Om))$,
  and of $(\veps)_{\eps\in (0,1)}$ and $(\pe(\veps))_{\eps\in (0,1)}$ in $L^\infty(\Om\times (0,\infty))$, from (\ref{44.3})
  and (\ref{44.4}) we immediately obtain (\ref{44.1}) and (\ref{44.2}) with some suitably large $C>0$.
\qed
As a preparation for our argument related to the continuity features claimed in Theorem \ref{theo10},
let us briefly record the following density property of the space appearing in the boundedness statements
from Lemma \ref{lem44}.
\begin{lem}\label{lem45}
  The set $W_N^{2,\infty}(\Om)$ is dense in $C^0(\bom)$.
\end{lem}
\proof
  This immediately follows from standard parabolic theory, which namely asserts that if 
%given any $\vp\in C^0(\bom)$
%  we can find a classical solution $z\in C^0(\bom\times [0,\infty)) \cap C^{2,1}(\bom\times (0,\infty))$ of
%  $z_t=\Del z$ in $\Om\times (0,\infty)$ with $\frac{\pa z}{\pa\nu}=0$ on $\pO$ and $z|_{t=0}=\vp$
  we let $(e^{t\Del})_{t\ge 0}$ denote
  the Neumann heat semigroup on $\Om$, then given any $\vp\in C^0(\bom)$ we have $e^{t\Del} \vp \in W_N^{2,\infty}(\Om)$
  for all $t>0$ and $e^{t\Del} \vp \to \vp$ in $C^0(\bom)$ as $t\searrow 0$.
\qed
As a consequence of the above a priori estimates, based on a standard
extraction procedure we can now construct a global solution in the sense of (\ref{dw}).
\begin{lem}\label{lem5}
  Let $n\ge 1$, and let (\ref{ie}) and (\ref{pe}) hold.
  Then there exist $(\eps_j)_{j\in\N} \subset (0,1)$ as well as nonnegative functions
  $u$ and $v$ on $\Om\times (0,\infty)$ 
  such that $\eps_j\searrow 0$ as $j\to\infty$, that (\ref{reg}) holds,
  and that as $\eps=\eps_j\searrow 0$ we have
  \begin{eqnarray}
	& & \ueps \wto u
	\qquad \mbox{in $L^2_{loc}(\bom\times (0,\infty))$,}
	\label{5.1} \\
	& & \veps \to v \mbox{ and } \na\veps\to\na v
	\qquad \mbox{a.e.~in } \Omega\times (0,\infty),
	\label{5.01} \\
	& & \veps\to v 
	\qquad \mbox{in } L^2_{loc}((0,\infty);W^{1,q}(\Om))
	\quad \mbox{for all $q\in [1,\frac{2n}{(n-2)_+})$,}
	\label{5.2} \\
	& & \veps(\cdot,t) \to v(\cdot,t)
	\qquad \mbox{in $W^{1,q}(\Om)$ for a.e.~$t>0$}
	\quad \mbox{for each $q\in [1,\frac{2n}{(n-2)_+})$,}
	\label{5.02} \\
	& & \na\veps \wsto \na v
	\qquad \mbox{in } L^\infty_{loc}((0,\infty);L^2(\Om)) 
	\qquad \qquad \mbox{and} 
	\label{5.22} \\
	& & \frac{\ueps\veps}{1+\eps\ueps} \wto uv
	\qquad \mbox{in } L^1_{loc}(\bom\times (0,\infty)).
	\label{5.23}
  \end{eqnarray}
  In the sense of Definition \ref{dw}, $(u,v)$ forms a global very weak solution of (\ref{0}) which satisfies
  \be{massu}
	\io u(\cdot,t) = \ouz |\Om|
	\qquad \mbox{for a.e. } t>0.
  \ee
  Moreover, for all $\tau>0$ there exists $C(\tau)>0$ such that
  \be{5.3}
	\io \Big| A^{-\frac{1}{2}}\big(u(\cdot,t)-\ouz\big)\Big|^2
	\le C(\tau)
	\qquad \mbox{for a.e. } t>\tau
  \ee
  and
  \be{5.33}
	\int_t^{t+1} \io u^2 
	+ \int_t^{t+1} \io |\Del v|^2
	+ \int_t^{t+1} \io |\na v|^4
	\le C(\tau)
	\qquad \mbox{for all } t>\tau,
  \ee
  and we have
  \be{5.4}
	A^{-\frac{1}{2}}(\ueps-\ouz) \wsto A^{-\frac{1}{2}}(u-\ouz)
	\qquad \mbox{in } L^\infty_{loc}((0,\infty);L^2(\Om))
  \ee
  as $\eps=\eps_j\searrow 0$.
\end{lem}
\proof
  From (\ref{4.2}) it follows that
  \bas
	(\ueps)_{\eps\in (0,1)}
	\mbox{ is bounded in $L^2(\Om\times (\tau,T))$ for all $\tau>0$ and $T>\tau$,}
  \eas
  while (\ref{4.3}), (\ref{4.4}), (\ref{4.01}) and (\ref{44.2}) guarantee that
  \bas
	& & (\veps)_{\eps\in (0,1)}
	\mbox{ is bounded in $L^2((\tau,T);W^{2,2}(\Om))$, in $L^4((\tau,T);W^{1,4}(\Om))$} \\[2mm]
	& & \hs{20mm}
	\mbox{and in $L^\infty((\tau,T);W^{1,2}(\Om))$ for all $\tau>0$ and $T>\tau$,}
  \eas
  and that
  \be{5.44}
	(v_{\eps t})_{\eps\in (0,1)}
	\mbox{ is bounded in } L^\infty\big((0,\infty); (W_N^{2,\infty}(\Om))^\star\big).
  \ee
  A standard extraction argument based on an Aubin-Lions lemma, and
  relying on the compactness of the embedding $W^{2,2}(\Om) \hra W^{1,q}(\Om)$ for 
  all $q\in [1,\frac{2n}{(n-2)_+})$, thus provides $(\eps_j)_{j\in\N}\subset (0,1)$ as well as nonnegative functions
  $u\in L^2_{loc}(\bom\times (0,\infty))$ and 
  $v\in L^2_{loc}((0,\infty);W^{2,2}(\Om))\cap L^4_{loc}((0,\infty));W^{1,4}(\Om)) \cap L^\infty_{loc}((0,\infty);W^{1,2}(\Om))$
  such that $\eps_j\searrow 0$ as $j\to\infty$, and that (\ref{5.1}), (\ref{5.01}), (\ref{5.2}), (\ref{5.02}) and (\ref{5.22})
  hold as $\eps=\eps_j\searrow 0$, where recalling (\ref{mass}), (\ref{vinfty}), (\ref{ie}), (\ref{4.1}) and (\ref{4.01}),
  and again using (\ref{4.2})-(\ref{4.4}), we readily obtain that also
  \be{5.45}
	u\in L^\infty((0,\infty);L^1(\Om))
	\qquad \mbox{and} \qquad
	v\in L^\infty(\Om\times (0,\infty)),
  \ee
  and that (\ref{5.33}), (\ref{5.3}) and (\ref{5.4}) as well as (\ref{massu}) hold.
  Furthermore, since $0\le \frac{\xi}{1+\eps\xi} \nearrow \xi$ as $\eps\searrow 0$ for all $\xi\ge 0$, and since thus the
  $L^1$ convergence feature trivially contained in (\ref{4.2}) ensures that also $\frac{\ueps}{1+\eps\ueps} \wto u$
  in $L^1_{loc}(\bom\times (0,\infty))$ as $\eps=\eps_j\searrow 0$ thanks to Lemma \ref{lem77},
  it also follows that (\ref{5.23}) holds,
  because $(\frac{\ueps\veps}{1+\eps\ueps})_{\eps\in (0,1)}$ is bounded in $L^2(\Om\times (\tau,T))$ and hence relatively
  compact in with respect to the weak topology in $L^1(\Om\times (\tau,T))$ for all $\tau>0$ and $T>\tau$ due to
  (\ref{4.2}) and (\ref{vinfty}), and because whenever $(\eps_{j_k})_{k\in\N}$ is a subsequence of $(\eps_j)_{j\in\N}$
  such that $\frac{\ueps\veps}{1+\eps\ueps} \wto z$ in $L^1_{loc}(\bom\times (0,\infty))$ 
  with some $z\in L^1_{loc}(\bom\times (0,\infty))$ as $\eps=\eps_{j_k} \searrow 0$, due to the pointwise approximation 
  property in (\ref{5.2}) a well-known result (\cite[Lemma A.1]{zhigun_surulescu_uatay}) becomes applicable so as to identify
  $z=uv$.\abs
  To derive the identities in (\ref{wu}) and (\ref{wv}) from this, we only need to observe that 
  for each $\vp\in C_0^\infty(\bom\times (0,\infty))$
  fulfilling $\frac{\pa\vp}{\pa\nu}=0$ on $\pO\times (0,\infty)$, according to (\ref{0eps}) we have
  \bas
	- \int_0^\infty \io \ueps \vp_t 
	= \int_0^\infty \io \ueps \pe(\veps) \Del\vp
	\qquad \mbox{for all } \eps\in (0,1)
  \eas
  and  
  \bas
	- \int_0^\infty \io \veps \vp_t 
	= \int_0^\infty \io \veps \Del\vp
	- \int_0^\infty \io \frac{\ueps\veps}{1+\eps\ueps} \vp
	\qquad \mbox{for all } \eps\in (0,1),
  \eas
  that (\ref{5.1}) and (\ref{5.2}) clearly entail that
  \bas
	\int_0^\infty \io \ueps \vp_t \to \int_0^\infty \io u\vp_t,
	\quad 
	\int_0^\infty \io \veps\vp_t \to \int_0^\infty \io v \vp_t
	\quad \mbox{and} \quad
	\int_0^\infty \io \veps \Del\vp \to \int_0^\infty \io v\Del \vp
  \eas
  as $\eps=\eps_j\searrow 0$, that due to (\ref{5.23}) we have
  \bas
	\int_0^\infty \io \frac{\ueps\veps}{1+\eps\ueps} \vp \to \int_0^\infty \io uv\vp
  \eas
  as $\eps=\eps_j\searrow 0$,
  and that (\ref{5.2}) together with (\ref{vinfty}), (\ref{ie}), (\ref{pe}), (\ref{phi}) and the dominated convergence 
  theorem ensure that $\pe(\veps)\to\phi(v)$ in $L^2_{loc}(\bom\times (0,\infty))$ and hence, again by (\ref{5.1}),
  \bas
	\int_0^\infty \io \ueps\pe(\veps) \Del\vp \to \int_0^\infty \io u\phi(v)\Del\vp
  \eas
  as $\eps=\eps_j\searrow 0$.\abs
  It remains to note that in addition to (\ref{5.44}) we know from Lemma \ref{lem44} that also
  \bas
	(u_{\eps t})_{\eps\in (0,1)}
	\mbox{ is bounded in } L^\infty\big((0,\infty); (W_N^{2,\infty}(\Om))^\star\big),
  \eas
  so that since both $L^1(\Om)$ and $L^\infty(\Om)$ are compactly embedded into $(W_N^{2,\infty}(\Om))^\star$,
  in view of (\ref{mass}), (\ref{vinfty}) and (\ref{ie}) we may twice employ the Arzel\`a-Ascoli theorem to infer that
  \be{5.7}
	\ueps\to u
	\quad \mbox{and} \quad
	\veps\to v
	\qquad \mbox{in } C^0_{loc}\big([0,\infty);(W_N^{2,\infty}(\Om))^\star\big)
  \ee
  as $\eps=\eps_j\searrow 0$, and that
  \be{5.77}
	u(\cdot,t) \to u_0
	\quad \mbox{and} \quad
	v(\cdot,t) \to v_0
	\quad \mbox{in } (W_N^{2,\infty}(\Om))^\star
	\qquad \mbox{as } t\searrow 0.
  \ee
  Indeed, since $W_N^{2,\infty}(\Om)$ is dense in $C^0(\bom)$ by Lemma \ref{lem45}, and since the inclusion
  $C_0^\infty(\Om) \subset W_N^{2,\infty}(\Om)$ entails density of $W_N^{2,\infty}(\Om)$ also in $L^1(\Om)$,
  from (\ref{5.7}), (\ref{5.45}) and (\ref{5.77}) it follows by means of a standard argument that actually
  $u\in C^0_{w-\star}([0,\infty);(C^0(\bom))^\star)$ and $v\in C^0_{w-\star}([0,\infty);L^\infty(\Om))$ with
  $u(\cdot,t) \wsto u_0$ in $(C^0(\bom))^\star$ and $v(\cdot,t)\wsto v_0$ in $L^\infty(\Om)$ as $t\searrow 0$.
\qed
Our main result on global solvability in the considered general framework thus becomes an evident consequence:\abs
\proofc of Theorem \ref{theo10}. \quad
  Since due to (\ref{phi}) and (\ref{init}) we can clearly choose $(\pe)_{\eps\in (0,1)}$ as well as $(u_{0\eps})_{\eps\in (0,1)}$
  and $(v_{0\eps})_{\eps\in (0,1)}$ such that (\ref{pe}) and (\ref{ie}) hold,
  the statement actually is part of what has been asserted by Lemma \ref{lem5}.
\qed
\mysection{Large time relaxation. Proof of Theorem \ref{theo11}}\label{sect4}
This section is devoted to the investigation of the large time relaxation feature of (\ref{0}) described in Theorem \ref{theo11}.
For this purpose, throughout this section we shall assume that $\phi$ 
satisfies (\ref{phi1}), so that in addition to (\ref{pe}), $(\pe)_{\eps\in (0,1)}$ can be chosen in such a way that
\be{pe1}
	\sup_{\eps\in (0,1)} \|\pe'\|_{L^\infty((0,M))} <\infty
	\qquad \mbox{for all } M>0.
\ee
In fact, we shall see that under this assumption, unlike in the general setting from Section \ref{sect3} we can 
control the integral on the right-hand side of (\ref{3.1}) in terms of suitably decaying quantities, based on the following
observation.
\begin{lem}\label{lem6}
  Let $n\ge 1$, and assume that (\ref{ie}), (\ref{pe}) and (\ref{pe1}) hold. Then given any $q>n$, one can find
  $\Gamma_3(q)>0$ such that
  \be{6.1}
	\frac{d}{dt} \io \Big| A^{-\frac{1}{2}}(\ueps-\ouz)\Big|^2
	+ \frac{1}{\Gamma_3(q)} \io (\ueps-\ouz)^2
	\le \Gamma_3(q) \|\na\veps\|_{L^q(\Om)}^2
	\qquad \mbox{for all $t>0$ and } \eps\in (0,1).
  \ee
%  for all $t>0$ and $\eps\in (0,1)$.
%
\end{lem}
\proof
  We abbreviate $c_1:=\|v_0\|_{L^\infty(\Om)}+1$ and may then rely on (\ref{pe1}) to fix $c_1>0$ such that
  \be{6.2}
	\big|\pe'(\xi)\big| \le c_2
	\qquad \mbox{for all $\xi\in [0,c_1]$ and any } \eps\in (0,1).
  \ee
  We moreover employ a Morrey-type estimate to see that thanks to our assumption $q>n$ we can find $c_3=c_3(q)>0$ fulfilling
  \be{6.3}
	\big|\vp(x)-\vp(y)\big|
	\le c_3 \|\na\vp\|_{L^q(\Om)}
	\qquad \mbox{for all } \vp\in C^1(\bom) \mbox{{\ and each $x, y\in\Om$}}.
  \ee
  On the right-hand side of (\ref{3.1}), recalling (\ref{mass}), (\ref{vinfty}) and (\ref{ie}) we can therefore estimate the 
  integrand according to
  \bas
	& & \hs{-30mm}
	\Big| \ouz \pe\big(\veps(x,t)\big) - \ov{\ueps(\cdot,t) \pe\big(\veps(\cdot,t)\big)} \Big| \\
	&=& \Bigg| \frac{1}{|\Om|} \cdot \bigg\{ \io \ueps(y,t) dy\bigg\} \cdot \pe\big(\veps(x,t)\big)
	- \frac{1}{|\Om|} \io \ueps(y,t) \pe\big(\veps(y,t)\big) dy \bigg| \\
	&=& \frac{1}{|\Om|} \io \ueps(y,t) \cdot \Big| \pe\big(\veps(x,t)\big)-\pe\big(\veps(y,t)\big) \Big| dy \\
	&\le& \frac{1}{|\Om|} \cdot \bigg\{ \io \ueps(y,t) dy \bigg\} \cdot 
		\sup_{y\in\Om} \Big| \pe\big(\veps(x,t)\big)-\pe\big(\veps(y,t)\big) \Big| \\
	&\le& \ouz c_2 \cdot \sup_{y\in\Om}\big|\veps(x,t)-\veps(y,t)\big|\\
	&\le& \ouz c_2 c_3 \|\na\veps(\cdot,t)\|_{L^q(\Om)}
	\qquad \mbox{for all $x\in\Om, t>0$ and } \eps\in (0,1).
  \eas
  From (\ref{3.1}) we hence obtain (\ref{6.1}) if we let
  $\Gamma_3(q):=\max\{\Gamma_2 \, , \, \Gamma_2 \ouz^2 c_2^2 c_3^2 |\Om|\}$.
\qed
Indeed, in low-dimensional cases the right-hand side of (\ref{6.1}) is, up to an expression already
known to decay integrably fast in time, essentially dominated by the
dissipation rate encountered in Lemma \ref{lem2}. 
More precisely, taking suitable linear combinations leads to the following.
\begin{lem}\label{lem7}
  Let $n\le 3$, and assume (\ref{ie}), (\ref{pe}) and (\ref{pe1}).
  Then there exist $b>0$ and $\Gamma_4>0$ such that
  \be{F}
	\Feps(t):= \io \Big| A^{-\frac{1}{2}}\big(\ueps(\cdot,t)-\ouz\big) \Big|^2
	+ b\io \big|\na\veps(\cdot,t)\big|^2,
	\qquad t\ge 0, \ \eps\in (0,1),
  \ee
  satisfies
  \be{7.1}
	\Feps'(t) 
	+ \frac{1}{\Gamma_4} \io (\ueps-\ouz)^2
	+ \frac{1}{\Gamma_4} \io |\na\veps|^4
	\le \Gamma_4 \io |\na\veps|^2
	\qquad \mbox{for all $t>0$ and } \eps\in (0,1).
  \ee
%  for all $t>0$ and $\eps\in (0,1)$.
%
\end{lem}
\proof
  Using that $\max\{n,2\} < \frac{2n}{(n-2)_+}$ due to our assumption that $n\le 3$,
  we can pick $q>\max\{n,2\}$ such that $q<\frac{2n}{(n-2)_+}$. We then let $\Gamma_3=\Gamma_3(q)$ be as accordingly provided
  by Lemma \ref{lem6}, and taking $\Gamma_1$ from Lemma \ref{lem2} and letting 
  \be{7.2}
	b:=\frac{1}{2\Gamma_1 \Gamma_3},
  \ee
  we can draw on the compactness of the first among the two continuous embedings $W^{2,2}(\Om) \hra W^{1,q}(\Om) \hra W^{1,2}(\Om)$
  to infer from an associated Ehrling inequality and standard elliptic regularity theory that there exists $c_1>0$ fulfilling
  \be{7.3}
	\Gamma_3 \|\na\vp\|_{L^q(\Om)}^2
	\le \frac{b}{2} \|\Del\vp\|_{L^2(\Om)}^2 
	+ c_1 \|\na\vp\|_{L^2(\Om)}^2
	\qquad \mbox{for all $\vp\in C^2(\bom)$ such that $\frac{\pa\vp}{\pa\nu}=0$ on } \pO.
  \ee
  Then from (\ref{6.1}) we obtain that for all $t>0$ and $\eps\in (0,1)$,
  \bas
%	& & \hs{-30mm}
	\frac{d}{dt} \io \big| A^{-\frac{1}{2}}(\ueps-\ouz)\big|^2
	+ \frac{1}{\Gamma_3} \io (\ueps-\ouz)^2 	
	&\le& \Gamma_3 \|\na\veps\|_{L^q(\Om)}^2 \\
	&\le& \frac{b}{2} \io |\Del\veps|^2
	+ c_1 \io |\na\veps|^2,
%	\qquad \mbox{for all $t>0$ and } \eps\in (0,1),
  \eas
  which when added to (\ref{2.1}) shows that with $(\Feps)_{\eps\in (0,1)}$ as in (\ref{F}) we have
  \bas
	& & \hs{-45mm}
	\Feps'(t)
	+ \frac{1}{\Gamma_3} \io (\ueps-\ouz)^2
	+ \frac{b}{2} \io |\Del\veps|^2
	+ \frac{b}{\Gamma_1} \io |\na\veps|^4 \\
	&\le& \frac{b}{2} \io |\Del\veps|^2
	+ c_1 \io |\na\veps|^2 \\
	& & + b\Gamma_1 \io (\ueps-\ouz)^2
	\qquad \mbox{for all $t>0$ and } \eps\in (0,1).
  \eas
  In view of (\ref{7.2}), this yields (\ref{7.1}) if we choose
  $\Gamma_4:=\max\{2\Gamma_3 \, , \, \frac{\Gamma_1}{b} \, , \, c_1\}$.
\qed
According to the decay feature of $t\mapsto \io |\na\veps|^2$ included in (\ref{gradv}),
an analysis of the damped linear ODI in (\ref{7.1}) already entails the stabilization property of $u$ 
claimed in Theorem \ref{theo11}, and beyond this also provides some further information on decay of the signal gradient.
\begin{lem}\label{lem8}
  Let $n\le 3$, and suppose that (\ref{ie}), (\ref{pe}) and (\ref{pe1}) hold.
  Then there exists a null set $N_\star\subset (0,\infty)$ such that $A^{-\frac{1}{2}}\big(u(\cdot,t)-\ouz\big)\in L^2(\Om)$
  for all $t\in (0,\infty)\sm N_\star$ with
  \be{8.1}
	A^{-\frac{1}{2}}\big(u(\cdot,t)-\ouz\big)
	\to 0
	\quad \mbox{in } L^2(\Om)
	\qquad \mbox{as } (1,\infty)\sm N_\star \ni t\to\infty,
  \ee
  and one can find $C>0$ such that
  \be{8.2}
	\int_1^\infty \io |\na\veps|^4 \le C
	\qquad \mbox{for all } \eps\in (0,1).
  \ee
\end{lem}
\proof
  Since $A^{-\frac{1}{2}}$ is continuous on $L^2_\perp(\Om)$, we can fix $c_1>0$ in such a way that
  \bas
	\|A^{-\frac{1}{2}}\vp\|_{L^2(\Om)}^2 \le c_1\|\vp\|_{L^2(\Om)}^2
	\qquad \mbox{for all } \vp\in L^2_\perp(\Om),
  \eas
  whence if we take $\Gamma_4$ from Lemma \ref{lem7} and let $c_2:=\frac{1}{c_1\Gamma_4}$, then from (\ref{7.1})
  we infer that for the functions in (\ref{F}) we have
  \be{8.3}
	\Feps'(t) + c_2 \Feps(t) + \frac{1}{\Gamma_4} \io |\na\veps|^4
	\le \heps(t):=(\Gamma_4 + bc_2) \io |\na\veps|^2
	\qquad \mbox{for all $t>0$ and } \eps\in (0,1).
  \ee
  If here we first neglect the nonnegative third summand on the left, then by means of a comparison argument we obtain that
  \bea{8.4}
	\Feps(t)
	&\le& \Feps(1) e^{-c_2(t-1)}
	+ \int_1^t e^{-c_1(t-s)} \heps(s) ds \nn\\
	&\le& c_3 e^{-c_1(t-1)} + \int_1^t e^{-c_1(t-s)} \heps(s) ds
	\qquad \mbox{for all $t>1$ and } \eps\in (0,1),
  \eea
  with $c_3:=\sup_{\eps\in (0,1)} \Feps(1)$ being finite due to (\ref{4.1}) and (\ref{4.01}).\abs
  In order to make appropriate use of this in the framework of the sparse topological information on the approximation 
  properties of $(\veps)_{\eps\in (0,1)}$, and especially of $(\ueps)_{\eps\in (0,1)}$, provided by Lemma \ref{lem5},
  we note that (\ref{4.01}) and (\ref{4.1}) particularly ensure that for all $t_0>1$, with $(\eps_j)_{j\in\N}$ as found there 
  we have $\na\veps\wsto \na v$ and 
  $A^{-\frac{1}{2}}(\ueps-\ouz) \wsto A^{-\frac{1}{2}} (u-\ouz)$ in $L^\infty((t_0,t_0+1);L^2(\Om))$ 
  as $\eps=\eps_j\searrow 0$, and that thus, according to lower semicontinuity of the norms in these spaces with respect to
  the considered convergence type, for
  \be{8.44}
	\F(t):= \io \Big|A^{-\frac{1}{2}}\big(u(\cdot,t)-\ouz\big) \Big|^2
	+ b \io \big|\na v(\cdot,t)\big|^2,
	\qquad t>0,
  \ee
  we have 
  \be{8.5}
	\|\F\|_{L^\infty((t_0,t_0+1))} \le \liminf_{\eps=\eps_j\searrow 0} \|\Feps\|_{L^\infty((t_0,t_0+1))}
	\qquad \mbox{for all } t_0>1.
  \ee
  Here the right-hand side can be controlled by combining (\ref{8.4}) with (\ref{5.2}) and (\ref{gradv}):
  Indeed, for $t_0>1$ and $t\in (t_0,t_0+1)$, in (\ref{8.4}) we can estimate
  \bas
	\int_1^t e^{-c_1(t-s)} \heps(s) ds
	= e^{-c_1 t} \int_1^t e^{c_1 s} \heps(s) ds
	\le e^{-c_1 t_0} \int_1^{t_0+1} e^{c_1 s} \heps(s) ds
	\qquad \mbox{for all } \eps\in (0,1),
  \eas
  where we may use that (\ref{5.2}) warrants that
  \be{8.6}
	\heps\to h
	\quad \mbox{in } L^1_{loc}((0,\infty))
	\qquad \mbox{as } \eps=\eps_j \searrow 0
  \ee
  with $h(t):=(\Gamma_4+bc_2) \io \big|\na v(\cdot,t)\big|^2$, $t>0$.
  Therefore, (\ref{8.5}) along with (\ref{8.4}) shows that
  \bas
	\|\F\|_{L^\infty((t_0,t_0+1))}
	&\le& c_3 e^{-c_2(t_0-1)} + e^{-c_1 t_0} \int_1^{t_0+1} e^{c_1 s} h(s) ds \\
	&=& c_3 e^{-c_2(t_0-1)} + e^{c_1} \int_1^{t_0+1} e^{-c_1(t_0+1-s)} h(s) ds
	\qquad \mbox{for all } t_0>1,
  \eas
  so that since $\int_1^\infty h(s) ds$ is finite by (\ref{gradv}) and (\ref{8.6}), and since thus
  \bas
	\int_1^{t_0+1} e^{-c_1(t_0+1-s)} h(s) ds 
	= \int_1^\infty \one_{(1,t_0+1)}(s) e^{-c_1(t_0+1-s)} h(s) ds 
	\to 0
	\qquad \mbox{as } t_0\to\infty
  \eas
  according to the dominated convergence theorem, it follows that
  \bas
	\|\F\|_{L^\infty((t_0,t_0+1))} \to 0
	\qquad \mbox{as } t_0\to\infty.
  \eas
  With some suitably chosen null set $N_\star\subset (0,\infty)$, this readily establishes the claimed conclusion together
  with (\ref{8.1}), whereas (\ref{8.2}) can be seen by going back to (\ref{8.3}) and integrating over $t\in (1,T)$ for $T>1$,
  which namely reveals that in view of our definition of $c_3$,
  \bas
	\frac{1}{\Gamma_4} \int_1^T \io |\na\veps|^4
	&\le& \Feps(1) + (\Gamma_4 + bc_2) \int_1^T \io |\na\veps|^2 \\
	&\le& c_3 + (\Gamma_4+bc_2) \cdot \frac{1}{2} |\Om| \cdot \big(\|v_0\|_{L^\infty(\Om)}+1\big)^2
	\qquad \mbox{for all $T>1$ and } \eps\in (0,1),
  \eas
  again thanks to (\ref{gradv}).
\qed
Finally, we only need 	%to appropriately combine the decay feature in (\ref{uv}) with (\ref{mass})	% and (\ref{8.2})
to observe that thanks to the fact that the exponent $4$ appearing in (\ref{8.2}) exceeds the currently considered spatial dimension, through a Morrey-type inequality the latter can be combined with (\ref{uv}) and (\ref{mass}) so as to 
yield decay of $v$ in the intended flavor.
\begin{lem}\label{lem9}
  If $n\le 3$ and (\ref{ie}), (\ref{pe}) as well as (\ref{pe1}) hold,
  then there exists a null set $N_{\star\star}\subset (0,\infty)$ such that 
  \be{9.1}
	v(\cdot,t)\to 0
	\quad \mbox{in } L^\infty(\Om)
	\qquad \mbox{as } (1,\infty)\sm N_{\star\star} \ni t\to\infty.
  \ee
\end{lem}
\proof
  Once more making explicit use of our restriction on $n$, we again employ a Morrey estimate to find $c_1>0$ such that
  \be{9.2}
	\Big\| \vp-\|\vp\|_{L^\infty(\Om)} \Big\|_{L^\infty(\Om)}
	\le c_1\|\na\vp\|_{L^4(\Om)}
	\qquad \mbox{for all } \vp\in W^{1,4}(\Om).
  \ee
  Moreover, combining (\ref{5.23}) with (\ref{uv}) and (\ref{5.01}) with (\ref{8.2}) we obtain that
  \bas
	\int_1^\infty \io uv<\infty
	\qquad \mbox{and} \qquad
	\int_1^\infty \io |\na v|^4 <\infty,
  \eas
  whence given $\eta>0$ we can choose $t_\eta>2$ suitably large fulfilling
  \be{9.3}
	\int_{t_\eta-1}^{t_\eta} \io uv \le \frac{m\eta}{2}
	\qquad \mbox{and} \qquad
	\int_{t_\eta-1}^{t_\eta} \big\|\na v(\cdot,t)\big\|_{L^4(\Om)} dt \le \frac{\eta}{2c_1},
  \ee
  where $m:=|\Om| \ouz$ is positive according to (\ref{init}).
  To see that, in fact, with some null set $N_{\star\star}\subset (0,\infty)$ we have
  \be{9.4}
	\|v(\cdot,t)\|_{L^\infty(\Om)} \le \eta
	\qquad \mbox{for all } t\in (t_\eta,\infty)\sm N_{\star\star},
  \ee
  we note that in view of the continuity of the embedding $W^{1,4}(\Om)\hra L^\infty(\Om)$ we know from (\ref{5.02}) that
  $\veps(\cdot,t)\to v(\cdot,t)$ in $L^\infty(\Om)$ for a.e.~$t>0$ as $\eps=\eps_j\searrow 0$, with $(\eps_j)_{j\in\N}$
  as provided there. As a consequence of this, namely, we may rely on (\ref{vinfty}) to see that with some
  null set $N_{\star\star}\subset (0,\infty)$,
  \be{9.5}
	\|v(\cdot,t)\|_{L^\infty(\Om)}
	\le \|v(\cdot,t_0)\|_{L^\infty(\Om)}
	\qquad \mbox{for all $t_0\in (0,\infty)\sm N_{\star\star}$ and each } t\in (t_0,\infty)\sm N_{\star\star},
  \ee
  while thanks to (\ref{massu}), upon enlarging $N_{\star\star}$ if necessary we may assume that moreover
  \be{9.6}
	\io u(\cdot,t) = m
	\qquad \mbox{for all } t\in (0,\infty)\sm N_{\star\star}.
  \ee
  In particular, using (\ref{9.2}) and (\ref{9.3}) we can estimate
  \bas
	& & \hs{-25mm}
	\int_{t_\eta-1}^{t_\eta} \|v(\cdot,t)\|_{L^\infty(\Om)} \cdot m dt \\
	&=& \int_{t_\eta-1}^{t_\eta} \io u(x,t) v(x,t) dxdt
	- \int_{t_\eta-1}^{t_\eta} \io u(x,t) \cdot \Big\{ v(x,t)-\|v(\cdot,t)\|_{L^\infty(\Om)}\Big\} dxdt \\
	&\le& \frac{m\eta}{2}
	+ \int_{t_\eta-1}^{t_\eta} \bigg\{ \io u(x,t) dx \bigg\} \cdot 
		\Big\| v(\cdot,t)-\|v(\cdot,t)\|_{L^\infty(\Om)} \Big\|_{L^\infty(\Om)} dt \\
	&\le& \frac{m\eta}{2}
	+ c_1 m \int_{t_\eta-1}^{t_\eta} \big\|\na v(\cdot,t)\big\|_{L^4(\Om)} dt \\
	&\le& \frac{m\eta}{2}
	+ c_1 m \cdot \frac{\eta}{2c_1},
  \eas
  so that since, on the other hand,
  \bas
	\int_{t_\eta-1}^{t_\eta} \|v(\cdot,t)\|_{L^\infty(\Om)} \cdot m dt
	\ge m \cdot {\rm{ess}} \hs{-3mm} \inf_{\hs{-5mm} t\in (t_\eta-1,t_\eta)} \|v(\cdot,t)\|_{L^\infty(\Om)}
  \eas
  by (\ref{9.5}) and (\ref{9.6}), it follows that
  \bas
	{\rm{ess}} \hs{-3mm} \inf_{\hs{-5mm} t\in (t_\eta-1,t_\eta)} \|v(\cdot,t)\|_{L^\infty(\Om)} \le \eta.
  \eas
  Agan thanks to (\ref{9.5}), this entails (\ref{9.4}) and thereby implies (\ref{9.1}).
\qed
Our main result on large time stabilization has thereby been achieved already.\abs
\proofc of Theorem \ref{theo11}. \quad
  Noting that our assumptions in (\ref{phi1}) enable us to choose $(\phi_\eps)_{\eps\in (0,1)}$ in such a way that both
  (\ref{pe}) and (\ref{pe1}) hold, we may take $(u,v)$ as accordingly provided by Lemma \ref{lem5}, and then conclude as intended by
  combining Lemma \ref{lem8} with Lemma \ref{lem9} and letting $N:=N_\star\cup N_{\star\star}$, with the null sets
  $N_\star$ and $N_{\star\star}$ introduced there.
\qed
\mysection{Appendix}
As we could not find an appropriate reference for this in the literature, let us finally include a derivation of the following 
general statement on weak $L^1$ convergence that has been used in Lemma \ref{lem5}.
\begin{lem}\label{lem77}
  Let $N\ge 1$ and $G\subset\R^N$ be measurable with $|G|<\infty$, and suppose that $(\rho_j)_{j\in\N} \subset C^0([0,\infty))$
  and $(w_j)_{j\in\N} \subset L^1(G;[0,\infty))$ are such that as $j\to\infty$ we have
  \be{77.1}
	\sup_{\xi\in [0,M]} \big|\rho_j(\xi)-\xi\big| \to 0
	\qquad \mbox{for all } M>0
  \ee
  and
  \be{77.2}
	w_j \wto w
	\qquad \mbox{in } L^1(G),
  \ee
  and that there exists $K>0$ such that
  \be{77.3}
	\big|\rho_j(\xi)\big| \le K\xi
	\qquad \mbox{for all $\xi\ge 1$ and } j\in\N.
  \ee
  Then 
  \be{77.4}
	\rho_j(w_j) \wto w
	\quad \mbox{in } L^1(G)
	\qquad \mbox{as } j\to\infty.
  \ee
  In particular, this conlcusion holds whenever $(\rho_j)_{j\in\N} \subset C^0([0,\infty))$ is such that
  \be{77.44}
	0 \le \rho_j(\xi) \nearrow \xi
	\quad \mbox{as $j\to\infty$ \qquad for all } \xi\ge 0.
  \ee
\end{lem}
\proof
  Since $(w_j)_{j\in\N}$ is relatively compact with respect to the weak topology in $L^1(G)$ by (\ref{77.2}), from the
  De la Vall\'ee-Poussin theorem we obtain $c_1>0$ and a function $\psi: [0,\infty)\to (0,\infty)$ such that 
  $\frac{\psi(\xi)}{\xi} \to +\infty$ as $\xi\to\infty$, and that
  \be{77.5}
	\int_G \psi(w_j) \le c_1
	\qquad \mbox{for all } j\in\N.
  \ee
  Now given $0\not\equiv \vp\in L^\infty(G)$ and $\eta>0$, we abbreviate $c_2:=\|\vp\|_{L^\infty(G)}$ and first pick
  $\del>0$ small enough fulfilling
  \be{77.6}
	\del\le \frac{\eta}{4c_1 c_2\cdot (K+1)},
  \ee
  then choose $M\ge 1$ such that
  \be{77.7}
	\frac{\psi(\xi)}{\xi} \ge \frac{1}{\del}
	\qquad \mbox{for all } \xi>M,
  \ee
  and taking $c_3>0$ large enough such that in accordance with (\ref{77.2}) we have
  \be{77.8}
	\int_G w_j \le c_3
	\qquad \mbox{for all } j\in\N,
  \ee
  we use (\ref{77.1}) and again (\ref{77.2}) to fix $j_\eta\in\N$ in such a way that
  \be{77.9}
	\big|\rho_j(\xi)-\xi\big|
	\le \frac{\eta}{4c_2 |G|}
	\qquad \mbox{for all $j\ge j_\eta$ and } \xi\in [0,M]
  \ee
  as well as
  \be{77.10}
	\bigg| \int_G w_j \vp - \int_G w\vp\bigg|
	\le \frac{\eta}{2}
	\qquad \mbox{for all } j\ge j_\eta.
  \ee
  Then for any $j\ge j_\eta$, in the identity
  \be{77.11}
	\int_G \rho_j(w_j)\vp - \int_G w\vp
	= \int_G \big\{ \rho_j(w_j)-w_j\big\}\cdot\vp
	+ \bigg\{ \int_G w_j \vp - \int_G w\vp\bigg\},
	\qquad j\in\N,
  \ee
  we can use our definition of $c_2$ to estimate
  \bea{77.12}
	\hs{-10mm}
	\bigg| \int_G \big\{ \rho_j(w_j)-w_j\big\}\cdot\vp \bigg|
	&\le& c_2 \int_G \big|\rho_j(w_j)-w_j\big| \nn\\
	&=& c_2 \int_{\{w_j\le M\}} \big|\rho_j(w_j)-w_j\big|
	+ c_2 \int_{\{w_j>M\}} \big|\rho_j(w_j)-w_j\big|
	\ \mbox{for all } j\in\N,
  \eea
  where according to (\ref{77.9}),
  \be{77.13}
	c_2 \int_{\{w_j\le M\}} \big|\rho_j(w_j)-w_j\big|
	\le c_2 \int_{\{w_j\le M\}} \frac{\eta}{4c_2 |G|}
	= \frac{\eta}{4} \cdot \frac{|\{w_j\le M\}|}{|G|} 
	\le \frac{\eta}{4}
	\qquad \mbox{for all } j\ge j_\eta.
  \ee
  Moreover, recalling that $M\ge 1$ we may rely on (\ref{77.3}) to see that thanks to (\ref{77.7}), (\ref{77.5}) and (\ref{77.6}),
  \bas
	c_2 \int_{\{w_j>M\}} \big|\rho_j(w_j)-w_j\big|
	&\le& c_2 \int_{\{w_j>M\}} \big|\rho_j(w_j)\big|
	+ c_2 \int_{\{w_j>M\}} w_j \\
	&\le& c_2 \cdot (K+1) \int_{\{w_j>M\}} w_j \\
	&=& c_2\cdot (K+1) \int_{\{w_j>M\}} \frac{w_j}{\psi(w_j)} \cdot \psi(w_j) \\
	&\le& c_2\cdot (K+1) \del \int_{\{w_j>M\}} \psi(w_j) \\
	&\le& c_2\cdot (K+1)\del \cdot c_1 \\
	&\le& \frac{\eta}{4}
	\qquad \mbox{for all } j\in\N.
  \eas
  Together with (\ref{77.13}) inserted into (\ref{77.12}), this shows that (\ref{77.11}) along with (\ref{77.10}) implies the
  inequality
  \bas
	\bigg| \int_G \rho_j(w_j) \vp - \int_G w\vp \bigg|
	\le \frac{\eta}{4} + \frac{\eta}{4}
	+ \bigg| \int_G w_j \vp - \int_G w\vp \bigg|
	\le \eta
	\qquad \mbox{for all } j\ge j_\eta,
  \eas
  from which (\ref{77.4}) follows for such $(\rho_j)_{j\in\N}$ due to the fact that $\eta>0$ and 
  $\vp\in L^\infty(G) \cong (L^1(G))^\star$ were arbitrary.\abs
  The additional claim concerning sequences $(\rho_j)_{j\in\N}$ fulfilling (\ref{77.44}) results from this upon observing that
  in this case (\ref{77.3}) is trivially satisfied, whereas (\ref{77.1}) is a consequence of Dini's theorem.
\qed
\bigskip

{\bf Acknowledgement.} \
  The first author was funded by the 
  China Scholarship Council (No. 202006630070). 
  The second author acknowledges support of the {\em Deutsche Forschungsgemeinschaft} (Project No.~462888149).
\small

\end{document}